\documentclass[twocolumn]{autart}
\usepackage{amssymb,amsfonts}
\usepackage[fleqn]{amsmath}
\usepackage{babel}
\usepackage{diymath}
\usepackage[pdftex]{hyperref}
\hypersetup{colorlinks=true,allcolors=auto,pdfborder={0 0 0}}
\edef\endfrontmatter{%
    \unexpanded\expandafter{\endfrontmatter}
    \noexpand\endNoHyper 
}

\journal{arXiv}
\usepackage{graphicx,subfigure}
\DeclareGraphicsExtensions{.pdf}
\graphicspath{{Figure/}}

\usepackage{savesym}
\savesymbol{AND}
\usepackage{algorithm}
\usepackage{algorithmic}

\usepackage{natbib}

\allowdisplaybreaks

\begin{document}
\begin{frontmatter}

	\title{Stabilization of continuous-time Markov/semi-Markov jump linear systems via finite data-rate feedback\thanksref{footnoteinfo}}

	\author[szucms]{Jingyi Wang}
    \ead{wangjingyi@szu.edu.cn}
	\author[szucms]{Jianwen Feng}
	\ead{fengjw@szu.edu.cn}
	\author[szucms]{Chen Xu}
    \ead{xuchen\_{szu}@szu.edu.cn}
	\author[whu]{Xiaoqun Wu}
    \ead{xqwu@whu.edu.cn}
	\author[bhu]{Jinhu L\"{u}}
    \ead{jhlu@iss.ac.cn}

	\address[szucms]{College of Mathematics and Statistics, Shenzhen University, Shenzhen 518060, PR China}
	\address[whu]{School of Mathematics and Statistics, Hubei Key Laboratory of Computational Science, Wuhan University, Wuhan, China}
	\address[bhu]{State Key Laboratory of Software Development Environment, School of Automation Science and Electrical Engineering, Beijing Advanced Innovation Center for Big Data and Brain Computing, Beihang University, Beijing, China}

\thanks[footnoteinfo]{This paper was not presented at any IFAC meeting. Corresponding author J.~Feng. Tel. +86-755-26532674.}

	\maketitle

	\begin{abstract}
		This paper investigates almost sure exponential stabilization of continuous-time Markov jump linear systems (MJLSs) under communication data-rate constraints by introducing sampling and quantization into the feedback control.
		Different from previous works, the sampling times and the jump times are independent \hl{of each other} in this paper.
		The quantization is recursively adjusted on the sampling time,
		and its updating strategy does not depend on the switching in a sampling interval.
		In other words, the explicit value of the switching signal in a sampling interval is not necessary.
		The numerically testable condition is developed to ensure almost sure exponential stabilization of MJLSs under the proposed communication and control protocols.
		We also drop the assumption of stabilizability of all individual modes  \hl{required  in previous works} about the switched systems.
		Moreover, we extend the result to the case of continuous-time semi-Markov jump linear systems (semi-MJLSs) via the semi-Markov kernel approach.
		Finally, some numerical examples are presented to illustrate the effectiveness on stabilization of the proposed communication and control protocols.
	\end{abstract}
\end{frontmatter}


%
%


\section{Introduction}\label{sec:introduction}



In recent decades, the switched system has attracted considerable attention because of its strong engineering background, see \cite{ieeecsm1999liberzon}, \cite{cdc1999hespanha:adt} and \cite{book:ssc}. 
In particular,
the Markov jump linear system has been considered as a special case of the switched system, when the switching signal is modeled as finite state Markov process,
because it has been widely used to model many practical systems with abrupt random changes, such as power systems \citep{auto2007li:mark}, freeway transportation system \citep{auto2017zhang:mjs}, networked control systems \citep{auto2013wang:mjs}, etc.




The stability and stabilization problems of MJLSs have received considerable attention in recent years ( e.g., \cite{auto2006bolzern,ieeetac2012li:eass,tac2016song:ass,tac2016zhanglx,acc2009zhang:mjls,auto2017gabriela,tac2018wu,ieeetac2018yang:ddrc} and the references therein).
\hl{The MJLS can be used to model practical systems with random changes.}
The continuous-time Markov process is used to describe the switching mechanism, including the switching patterns, the jump (switching) times, and the sojourn (holding) times of the active mode.
The stability analysis of continuous-time MJLSs was considered in the almost sure sense and a sufficient condition was derived depending on the transition matrix in \cite{auto2006bolzern}, while the discrete-time case was presented in \cite{iet2013song:ass}.
Moreover, the switching MJLS was studied and its almost sure exponential stability was obtained, in which the transition rate matrix for the random Markov process was varied when a deterministic switching occurred in \cite{tac2016song:ass}.
The stability and stabilization of discrete-time semi-MJLSs were considered via the semi-Markov kernel approach in \cite{tac2016zhanglx}.
The second order stabilization problems of MJLSs were studied via explicitly constructing the stabilizing logarithmic quantizer and controller in \cite{acc2009zhang:mjls}.
The stability problem of semi-Markov and Markov switched systems was investigated by using the probability analysis method in \cite{tac2018wu}.
However, it is often a key restriction that the jump times of Markov process and the sampling times of the sampler are the same time series, which leads to the fact that the sampling and switching are simultaneous, but the switching is usually stochastic.
When the sampler don't know the switches occurring, can the system achieve stabilization via the controller?



%

Understanding control over communication networks was listed as a major challenge for the controls field \citep{csm2003murray}.
In engineering systems, the total communication capacity in bits per second may be large in the overall system, but each component is effectively allocated only a small portion \citep{ieeeprog2007nair:drc}.
The finite data-rate feedback here means that measurement information transferred though a communication channel with finite bandwidth from the sensor to the controller (see Fig. \ref{fig:infostru}).
The finite data-rate feedback\hl{,  which can balance the communication capacity and control performance,} combines the reliable transmission of information in communication theory and feedback control of information in control theory.
Sampling and quantization are fundamental tools to deal with finite data-rate feedback problems in the modern control systems.
Actually, sampling is the reduction of a continuous-time signal to a discrete-time signal at the sampling time, and quantization is a kind of mapping from continuous signals to discrete sets by the prescribed rules.
%
%
With this motivation, stabilization of the control systems via finite data-rate feedback control was studied in {the continuous-time \citep{tac2000brockett:qfs,ieeetac2003liberzon:slsli,tac2005libhes,csl2021berger}, discrete-time \citep{ieeetac2001elia:slsli,auto2003liberzon,tcyb2019zhang} and switching \citep{ieeetac2017wakaiki:swls} settings, which may be subject to external disturbances \citep{ieeetac2012sharon:cq} or in nonlinear systems \citep{fss2019zhang,ijbnc2017shi}.}
%
Moreover, the techniques can be used to deal with the stabilization of systems with additive Gaussian white noise \citep{ieeetac2003liberzon:slsli}.

To the best of our knowledge, few studies have been conducted on the stabilization problem of Markov jump linear systems under finite data-rate feedback. 
In this paper, we consider the stabilization problem of the Markov/semi-Markov jump linear systems under quantized state feedback subject to communication data-rate constraints.
The main contributions of this paper are listed as follows.
\hl{First, we give the method to design the communication and control protocols and update the parameters of the protocol under communication data-rate constraints.
    Second, we derive testable sufficient conditions for the almost sure exponential stabilization of the MJLS under the proposed protocol, and drop the assumption of stabilizability of all individual modes.
    Third, the results are extended to the almost sure exponential stabilization of the semi-MJLSs.}

The structure of this article is listed as follows:
Section \ref{sec:preliminaries} introduces some preliminaries, including the model description of MJLSs, information patterns of the system, and some related concepts of almost sure exponential stabilization.
In Section \ref{sec:mr}, we obtain the sufficient conditions of almost sure exponential stabilization, which are dependent on the generator of the Markov process. Next, the updating rule of the quantization parameters is designed.
Moreover, the results are extended to the case of semi-MJLSs.
Some examples are given in Section \ref{sec:ns} to illustrate the effectiveness of our results.



%
%
%


\section{Preliminaries}\label{sec:preliminaries}

\subsection{Notations}

First, we write down the notations that will be used throughout this paper.
Let $\sR = (-\infty, +\infty)$ be the set of real numbers,
$\sR_{+} = [0, +\infty)$ be the set of non-negative real numbers,
$\mathbb{N} = \{0, 1, 2, 3, \ldots\}$ be the set of natural numbers,
$\mathbb{N}_{+} = \{1, 2, 3, \ldots\}$ be the positive integers,
$\mathbb{N}_{[m, n]} = \{m, m+1, m+2, \ldots, n\}$ where $ m<n $ for $m, n \in \sN$,
$\sR^n$ be the $n-$dimensional Euclidean space and $\sR^{n \times n}$ be the set of all $n \times n$ real matrices.
Let $\mathbf{1}_n=[1, 1, \ldots, 1]^{\top} \in \sR^n$ be the one vector, $\mI_n$ be the $n-$dimensional identity matrix,
$\mO$ be the zero matrix.
For $A,B \in \sR^{m \times n}$, let $[A,B]$ ($[A;B]\hl{=[A^{\top},B^{\top}]^{\top}}$) denote the horizontal (vertical) concatenation of $A$ and $B$.
Let the superscript $\top$  denote the transpose of a matrix. Let $\hl{\prod_{i=1}^{k}}B_i=B_k B_{k-1} \cdots B_1$ denote the left product of matrix $B_i$ ($i=1, 2, \ldots, k$). Let $\underline{\lambda}(\cdot)$ and $\bar{\lambda}(\cdot)$ denote the smallest and the largest eigenvalue of a symmetric matrix, respectively.
Let $\|\cdot\|$ denote $l_\infty$ norm, i.e., $\|x\|=\max _{1 \leq i \leq n}\left|x_{i}\right|$ on $\mathbb{R}^{n}$ and the corresponding induced matrix norm $\|A\|=\max _{1 \leq i \leq n} \sum_{j=1}^{m}\left|A_{i j}\right|$ on $\mathbb{R}^{n \times m}$. Let a triple $(\Omega, \mF , \pP )$ be the complete probability space, where
$\Omega$ represents the sample space,
$\mF$ is the $\sigma$-algebra of subsets of the sample space, known as the event space,
and $\pP$ is the probability measure on $\mF$ and the measure $\pP\{E\}$ is known as the probability of the event $E \in \mF$.

\subsection{Model description}\label{sec:md}

We will consider the stabilization problem for continuous-time	Markov jump linear systems as follows:
\begin{equation}\label{sys:init}
	\dot{x}(t)=A_{\sigma(t)}x(t) + B_{\sigma(t)}u(t), \quad x(0) = x_0,
\end{equation}
where $x(t) \in \sR^n$ is the state vector, $u(t) \in \sR^m$ is the control input, $A\in \sR^{n \times n}$, $B \in \sR^{n \times m}$ , and the switching signal $\sigma(t): \sR_{+} \rightarrow \sM$ is a finite-state, homogeneous Markov process taking value in a finite index set $\sM=\{1, 2, \ldots, M\}$ with   generator $\Gamma=[\gamma_{ij}] \in \sR^{M \times M}$ given by
\begin{equation*}
	\mprob{\sigma(t+h)=j | \sigma(t)=i}=
	\left\{\begin{aligned}
		 & \gamma_{ij} h+o(h),   &  & i \neq j, \\
		 & 1+\gamma_{ii} h+o(h), &  & i=j,
	\end{aligned}\right.
\end{equation*}
where $h>0$ and $\lim _{h \rightarrow 0} o(h) / h=0$, $\gamma_{ij} \geq 0$ ($i \neq j$) stands for
the transition rate from mode $i$ to mode $j$, and $\gamma_{ii}=-\sum_{j=1, j \neq i}^{M} \gamma_{ij}$
, which specifies the active mode $\sigma(t)$ at each time $t$.
The set of the matrix pair $\{(A_p, B_p), p \in \sM\}$ denotes a collection of matrices defining the modes.

Define the sequence of jump times of $\{\sigma(t), t \in \sR_{+}\}$ by
\begin{equation*}
	\tau_{0}=0,~  \tau_{k+1}=\inf \left\{t : t \geq \tau_{k}, \sigma(t) \neq \sigma({\tau_{k})}\right\},~  k \in \sN,
\end{equation*}
where we adopt the convention $\inf \emptyset=\infty$.
Define the sequence of sojourn times (holding times) of $\{\sigma(t), t \in \sR_{+}\}$ by
\begin{equation*}
	h_{k}= \tau_{k}-\tau_{k-1},~  k \in \sN_{+},
\end{equation*}
where $h_{k}=\infty$ if $\tau_{k}=\infty$.
Let $\sigma({\infty})=\sigma({\tau_{k-1}})$ if $\tau_{k}=\infty$.
The jump chain induced by  $\{\sigma(t), t \in \sR_{+}\}$ is defined to be
\begin{equation*}
	\sigma_{k}=\sigma(\tau_{k}),~  k \in \sN.
\end{equation*}
The sequence $\left\{\sigma_{k}, k \in \sN \right\}$ is called the embedded chain of $\{\sigma(t), t \in \sR_{+}\}$,
which is the discrete-time Markov process with transition probability matrix $\Lambda = [\lambda_{ij}] \in \sR^{M \times M}$ defined by $\lambda_{ij} = -\gamma_{ij}/\gamma_{ii}$ if $i \neq j$, $\lambda_{ii}=0$ otherwise. The sojourn times $h_{1}, h_{2}, h_{3},\ldots $ are independent exponential random variables with parameters $\gamma_{\sigma_{0}}, \gamma_{\sigma_{1}}, \gamma_{\sigma_{2}},\ldots, $ respectively, where $\gamma_{i}=-\gamma_{ii}$.
In other words, $\mprob{\tau_{k+1}-\tau_{k} < t | \sigma(\tau_{k})=i, \sigma(\tau_{k+1})=j}=1 - e^{-\gamma_{i} t}$ is independent of mode $j$ for any $t \in \sR_{+}$  (see, e.g., \cite{book:mp}).

\begin{assum}[Markov process]\label{asm:markproc}
	The Markov process $\{\sigma(t), t \in \sR_{+}\}$ is irreducible and aperiodic.
\end{assum}

Assumption \ref{asm:markproc} implies that the Markov process is ergodic and has a unique stationary distribution
$\pi =[\pi_1, \pi_2, \ldots, \pi_M]$
which can be calculated by
$\pi \Lambda = \pi$ and $\sum_{i=1}^M\pi_i=1$ (see, e.g., \cite{book:mp}).

The objective is to stabilize the Markov jump linear systems with the controller $u(t)$ under communication data-rate constraints.
In the sequel, the concept of stability is given as follows:

\begin{defn}\citep{book:sde}\label{def:as}
	\hl{The solution of system \eqref{sys:init} is said to be \textit{almost surely exponentially stable} if} there exists a positive constant $\varepsilon$ such that
	\begin{equation*}
		\mProb{\limsup_{t \rightarrow \infty} \sfrac{1}{t} \ln \norm{x(t)} \leq -\varepsilon} = 1,
	\end{equation*}
	for any  $x_0\in \sR^n$.
\end{defn}


Our another goal is to drop the assumption of \hl{stabilizability} of all individual modes, which is needed in the context of stabilization of finite data-rate feedback, e.g.,  \cite{auto2014liberzon:findatarate,ieeetac2018yang:ddrc}.
Let $\sM_{s}$ be the index set of the stabilizable pair $(A_p, B_p)$, i.e., there exists a state feedback gain matrix $K_p$ such that $A_p + B_p K_p$ is Hurwitz,
and Let $\sM_{u}$ be the index set of the unstabilizable pair $(A_p, B_p)$, i.e., there are not any matrices $K_p$ such that $A_p + B_p K_p$ is Hurwitz.
Obviously, $\sM=\sM_{s} \cup \sM_{u}$ and $\sM_{s} \cap \sM_{u} = \emptyset$.

If the pair $(A_p, B_p)$ is \hl{stabilizable}, we assume that the suitable stabilizing gain matrix $K_p$ has been selected and fixed such that $A_{p}+B_{p}K_{p}$ is Hurwitz, the matrix $K_{p}$ can be obtain, e.g., by solving some algebraic Riccati equations, and if the pair $(A_p, B_p)$ is \hl{unstabilizable}, we assume matrix $K_p = \mO$ with suitable dimension.

\subsection{Information Patterns}\label{sec:ip}

\begin{figure}[!htbp]
    \centering
    \includegraphics[width=\linewidth]{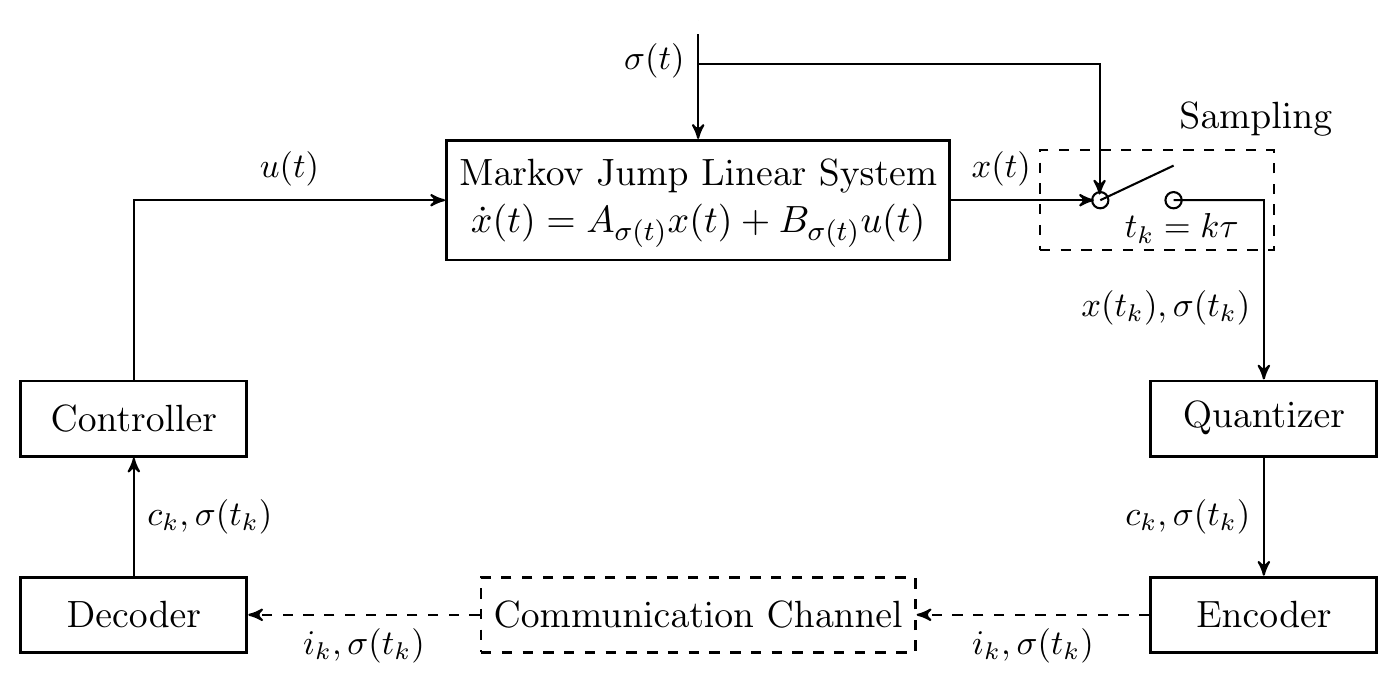}
    \caption{Block diagram of the MJLS and information pattern}
    \label{fig:infostru}
\end{figure}

In this paper, the controller is separated from the actuator and the sensor used to measure the system state,
and  the communication channel is noiseless.
The state information is processed and transferred similarly to \cite{auto2003liberzon,ieeetac2004mitter:ccc} and \cite{ieeetac2018yang:ddrc} as shown in Fig. \ref{fig:infostru} in the following standard way.
\begin{enumerate}
	\item Sampling: State measurements are taken at time $t_k=k\tau$, $k=0, 1, 2, \ldots$, where $\tau$ is a fixed sampling period.
	\item Quantizing and encoding: Each state measurement $x(t_k)$ is quantized and encoded by an integer $i_k$ from $1$ to $N^n$ by some rule given in below, where $N$ is an odd positive integer. In addition, the pair $(i_k, \sigma(t_k))$ \hl{is  encoded as a sequence of  bits, and sent to the decoder by the digital communication channel.}
	\item Decoding: The state $c_k$ and $\sigma(t_k)$ are decoded from \hl{the bitstream  of the pair $(i_k, \sigma(t_k))$ by the rules given in advance.}
	\item Controlling: The control signal $u(t)$ is then determined solely from the decoder’s state $(c_{k}, \sigma(t_{k}))$ according to the control protocol.
\end{enumerate}

\begin{figure}[!htbp]
	\centering
	\includegraphics[width=0.95\linewidth]{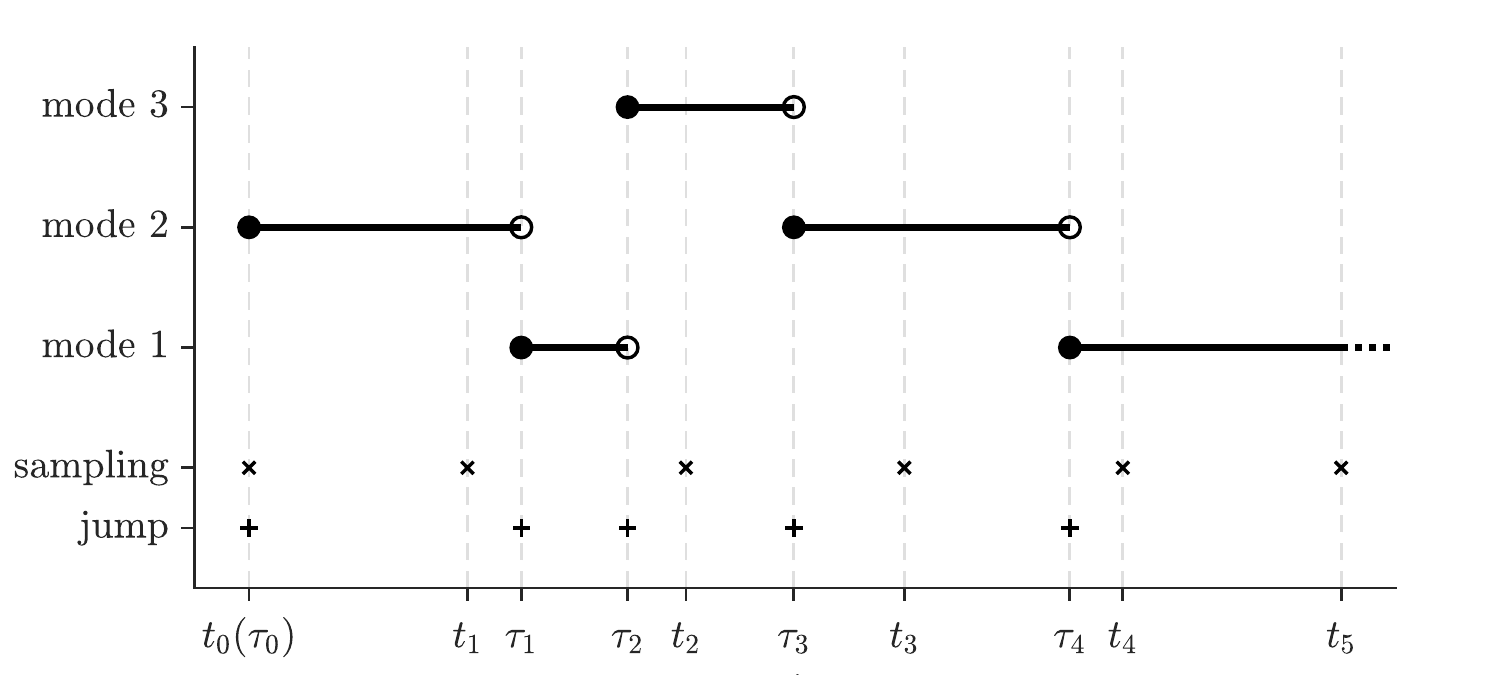}
	\caption{Sample path of a Markov process $\{\sigma(t), t \in \sR_{+}\}$ with  jump times and switching modes ($M=3$), and sampling times.}
	\label{fig:timemodel}
\end{figure}

\begin{rem}
	The sequences $\{\sigma(t_{k}), k \in \sN\}$ and $\{\sigma(\tau_{k}), k \in \sN\}$ are different. In particular, $\{\sigma(t_{k}), k \in \sN\}$ is dependent on the sampling time while $\{\sigma(\tau_{k}), k \in \sN\}$ is dependent on the switching. A connection between $\{\sigma(t_{k}), k \in \sN\}$ and $\{\sigma(\tau_{k}), k \in \sN\}$ is $\sigma(t_{k})=\sigma(\tau_{k'})$ where $\tau_{k'} = \max\{\tau_{i} : \tau_{i} \leq t_{k}, i \in \sN\}$.
	An illustration of such a switching pattern for the case $\sM=\{1,2,3\}$ is depicted in Fig. \ref{fig:timemodel}.
\end{rem}

Processing similarly to \cite{auto2003liberzon,ieeeprog2007nair:drc,ieeetac2004mitter:ccc,ieeetac2018yang:ddrc}, the data-rate (also known as bit-rate) between the encoder and the decoder
\begin{equation*}
	R = ({\log_{2}(N^{n}+1) + \log_{2}\left| \sM \right|})/{\tau}
\end{equation*}
is bits per unit of time, where $\left| \sM \right|$ is the cardinality of the index set $\sM$ (i.e., the number of modes).
$R$ is used to characterize the capacity of communications.

%
%

%

%
%



\begin{assum}[finite data-rate]\citep{auto2014liberzon:findatarate}\label{asm:sampletime}
	The sampling period $\tau$ satisfies
	$\Lambda_p < N$ for any $p \in \sM$, where $\Lambda_p = \norm{\exp(A_p\tau)}$.
\end{assum}

\begin{rem}
    \hl{The assumption is viewed as a constraint of data-rate,
    since the inequality $  \norm{\exp(A_p\tau)} < N $   requires $\tau$ to be sufficiently small
    with respect to $ N $.
    Combining the definition of $ R $, one can know that the data-rate $ R $ has a lower bound.  }
\end{rem}

\section{Main Results}\label{sec:mr}

The control objective is to stabilize the system defined in Section \ref{sec:md} in the sense of Definition \ref{def:as} 
while respecting the communication data-rate constraints described in Section \ref{sec:ip}.
Our results were inspired by the work of \cite{auto2014liberzon:findatarate} and \cite{ieeetac2018yang:ddrc}, where
all individual modes are stabilizable and switches actually occur less often than once per sampling period.

\subsection{Communication and control protocols}\label{sec:comcont}
%
%

In this subsection, we describe the communication and control strategy similarly to \cite{auto2014liberzon:findatarate} and \cite{ieeetac2018yang:ddrc}.


The initial state $x_{0}$ is unknown. At $t_{0}=0, $ the sensor and the controller are both provided with $x_{0}^{*}=0$ and arbitrarily selected initial estimates $E_{0}>0$ and $\delta_{0}>0$. Starting from $t_{0}=0$, at each sampling time $t_{k}$, the sensor determines if the state $x(t_{k})$ is inside the hypercube of radius $E_{k}$ centered at $x_{k}^{*}$ denoted by
\begin{equation}\label{set:encdec}
	\mS_k= \{x \in \sR^n : \norm{x-x^{\ast}_k} \leq E_k\}
\end{equation}
The hypercube $\mathcal{S}_{k}$ is the approximation of the reachable set at $t_{k}, $ which is also used as the range of quantization.
How to update $x^{\ast}_{k}$ and $E_{k}$ such that \eqref{set:encdec} holds i.e.,
\begin{equation}\label{eq:relxE}
	\left\|x(t_{k})-x_{k}^{*}\right\| \leq E_{k},
\end{equation}
will be given in Section \ref{sec:uq}.

At each sampling time $t_{k}$, the quantizer divides the hypercube $\mathcal{S}_{k}$ into $N^{n}$ equal hypercubic boxes, $N$ per dimension, each box is encoded by a unique integer index from $1$ to $N^{n}$, and the index $i_{k}$ of the box containing $x(t_{k})$ and the active mode $\sigma(t_{k})$ are transmitted to the decoder.
The decoder follows the same predefined indexing protocol as the encoder, so that it is able to reconstruct the center $c_{k}$ of the hypercubic box containing $x(t_{k})$ from $i_{k}$.
The controller then generates the control input
\begin{equation}\label{eq:continput}
	u(t) = K_{\sigma(t_{k})}\hat{x}(t)
\end{equation}
for $t \in\left[t_{k}, t_{k+1}\right)$, where $K_{\sigma(t_{k})}$ is the  feedback gain matrix, and $\hat{x}(t)$ is the state of the auxiliary system described by
\begin{equation}\label{eq:auxilinsys}
	\dot{\hat{x}}(t) = A_{\sigma(t_{k})}\hat{x}(t) + B_{\sigma(t_{k })} u(t)
\end{equation}
with the boundary condition $\hat{x}(t_k) = c_k$.
\hl{The auxiliary system is to design the feedback controller with the estimated state of $x(t)$, which is frequently unobservable.}
Simple calculation shows that
\begin{align}\label{eq:aux}
	\norm{{x}(t_k) - c_k} \leq \sfrac{1}{N}E_k \text{~and~}
	\norm{x^{\ast}_k - c_k} \leq \sfrac{N-1}{N}E_k.
\end{align}
\hl{Let} $\Delta_{k} = x^{\ast}_{k} - c_{k}$.

\subsection{Updating the parameters of the quantization} \label{sec:uq}
In the following, we will give the rules to update $x^{\ast}_{k}$ and $E_{k}$ such that equation \eqref{eq:relxE} holds.

Depending on whether switches occur, two case need to be considered respectively. Consider the sample interval $[t_k, t_{k+1})$, without loss of generality, \hl{let $\sigma(t_k)=p$ and $\sigma(t_{k+1})=q$  throughout this paper}.

\textit{Case 1: Sampling interval without switches occurring.}
\hl{Let} $e(t) = x(t) - \hat{x}(t)$.
Obviously,
$\dot{e}(t) = A_pe(t)$ and
$e(t_{k+1}^-) =\exp(A_p\tau)e(t_k) $
from equation \eqref{sys:init} and \eqref{eq:auxilinsys}.
One can obtain that
\begin{equation*}
	\norm{e(t_{k+1}^-)} = \norm{\exp(A_p\tau)e(t_k)} \leq \norm{\exp(A_p\tau)}\norm{e(t_k)}
\end{equation*}
So we can update $x^{\ast}_{k+1}$ and $E_{k+1}$ by using
\begin{subequations}\label{eq:updaxE1}
	\begin{align}
		 & x^{\ast}_{k+1} := \hat{x}(t_{k+1}^-) = \exp\big((A_p + B_p K_p) \tau\big)c_k, \label{eq:updaxE1a} \\
		 & E_{k+1} := \sfrac{\Lambda_p}{N} E_k. \label{eq:updaxE1b}
	\end{align}
\end{subequations}

\textit{Case 2: Sampling interval with switches occurring.}
Let $\tau_{k}^{1}$, $\tau_{k}^{2}$, $\ldots$, $\tau_{k}^{s}$
be the switching times of the Markov process, where $ s $ is the number of switches in $[t_{k}, t_{k+1})$.
Obviously, ${\tau}_{k}^{i}$ and $s$ are unknown.
%
Let $\tau_{k}^{0}=\tau_{k}$ and $\tau_{k}^{s+1}=\tau_{k+1}$.
From \eqref{sys:init} \eqref{eq:continput} and \eqref{eq:auxilinsys}, one can obtain that
\begin{equation}\label{equ:xtk1}
	x(t^{-}_{k+1}) = [\mI_{n}, \mO_{n}] \prod^{s}_{i=0}\exp\big(A_{p, \sigma(\tau_{k}^{i})}(\tau_{k}^{i+1}-\tau_{k}^{i})\big) [ x_{k}; c_{k}],
\end{equation}
where
\begin{equation*}
	A_{p, q} = \begin{bmatrix}
		A_q     & B_qK_p     \\
		\mO_{n} & A_p+B_pK_p \\
	\end{bmatrix}.
\end{equation*}

In order to estimate $x^{\ast}_{k+1}$ and $ E_{k+1} $, one can select some times $\tilde{\tau}_{k}^{1},  \tilde{\tau}_{k}^{2}, \ldots,  \tilde{\tau}_{k}^{M-1}$ as \textit{expected} switching times.
$\sigma(\tilde{\tau}_{k}^{i})$ is the \textit{expected} switching mode corresponding to the switching time $\tilde{\tau}_{k}^{i}$.
Let $\tilde{\tau}_{k}^{0}=\tau_{k}$ and $\tilde{\tau}_{k}^{M}=\tau_{k+1}$.
Processing similarly to \eqref{equ:xtk1}, one can obtain that
\begin{equation}\label{equ:xtk2}
	\tilde{x}(t^{-}_{k+1}) = [\mI_{n}, \mO_{n}] \prod^{M}_{i=0}\exp\big(A_{p, \sigma(\tilde{\tau}_{k}^{i})} (\tilde{\tau}_{k}^{i+1}-\tilde{\tau}_{k}^{i})\big) [c^{\top}_{k}, c^{\top}_{k}]^{\top},
\end{equation}
where set $\{\sigma(\tilde{\tau}_{k}^{0}),\ldots \sigma(\tilde{\tau}_{k}^{M-1})\} = \sM$ and  $\tilde{\tau}_{k}^{i+1} - \tilde{\tau}_{k}^{i} =\tau(1/\gamma_{p})/(\sum^m_{i=1}1/\gamma_{i})$ if $\sigma(\tilde{\tau}_{k}^{i})= p$.
Moreover, one can select some instants $ \check{\tau}_{k}^{1},\check{\tau}_{k}^{2}, \ldots,\check{\tau}_{k}^{w}$ as the \textit{worst} switching times.
Let $\check{\tau}_{k}^{0}=\tau_{k}$ and $\check{\tau}_{k}^{w+1}=\tau_{k+1}$.
Processing similarly to \eqref{equ:xtk1}, one can obtain that
\begin{equation}\label{equ:xtk3}
	\check{x}(t^{-}_{k+1}) = [\mI_{n}, \mO_{n}] \prod^{w}_{i=1}\exp\big(A_{p, \sigma(\check{\tau}_{k}^{i})}(\check{\tau}_{k}^{i+1}-\check{\tau}_{k}^{i})\big)
	[c^{\top}_{k}, c^{\top}_{k}]^{\top},
\end{equation}
where
\begin{multline*}
	\{w,\check{\tau}_{k}^{i},\sigma(\check{\tau}_{k}^{i}), i \in \sN_{[1, w]}\} \\
	=\argmax
	\big\|\prod^{w}_{i=1}\exp\big(A_{p, \sigma(\check{\tau}_{i})} (\check{\tau}_{k}^{i+1}-\check{\tau}_{k}^{i})\big) \big\|
\end{multline*}
can be seen as the outliers of the switching times.

Let $ \tilde{S}_{k} = [\mI_{n}, \mO_{n}] \prod^{M}_{i=0}\exp\big(A_{p, \sigma(\tilde{\tau}_{k}^{i})} (\tilde{\tau}_{k}^{i+1}-\tilde{\tau}_{k}^{i})\big)$
and $\check{S}_{k} = [\mI_{n}, \mO_{n}] \prod^{w}_{i=1}\exp\big(A_{p, \sigma(\check{\tau}_{k}^{i})} (\check{\tau}_{k}^{i+1}-\check{\tau}_{k}^{i})\big)  $.
From \eqref{eq:relxE} and \eqref{eq:aux}, one can get that
\begin{equation}\label{equ:equxhatx}
	\|x(t_{k+1}) - \check{x}(t^{-}_{k+1})\| \leq \| \check{S}_{k} \| \big(2\|x^{\ast}_{k}\|+ \sfrac{2N-1}{N}E_{k}\big).
\end{equation}
Moreover, from \eqref{eq:aux} and \eqref{equ:equxhatx}, by the triangle inequality one can get that
\begin{align}\label{eq:hold4}
	     & \|x(t_{k+1}) -\tilde{x}(t_{k+1})\| \notag                                                   \\
	\leq & (2 \| \check{S}_{k} \| + \|\tilde{S}_{k} - \check{S}_{k}\|) \|x^{\ast}_{k}\|\notag          \\
	     & + \sfrac{(N-1) \|\tilde{S}_{k} - \check{S}_{k}\| + (2N-1)\|\check{S}_{k}\|}{N}E_{k}.
\end{align}

So we can update $x^{\ast}_{k+1}$ and $E_{k+1}$ when switches occur by using
\begin{subequations}\label{eq:updaxE2}
	\begin{align}
		x^{\ast}_{k+1} := & \tilde{S}_{k} [c^{\top}_k, c^{\top}_{k}]^{\top}, \label{eq:updaxE21a} \\
		E_{k+1} :=        & \chi_k\|x^{\ast}_{k}\| + \psi_{k}E_{k},
	\end{align}
\end{subequations}
where
\begin{align*}
	\chi_k =  & 2 \| \check{S}_{k} \| + \|\tilde{S}_{k} - \check{S}_{k}\|,                 \\
	\psi_{k}= & ({(N-1) \|\tilde{S}_{k} - \check{S}_{k}\| + (2N-1)\|\check{S}_{k}\|})/{N}.
\end{align*}

\begin{rem}
    How to select $\{\check{\tau}_{k}^{i}, i \in \sN_{[1, M]}\}$ is challenging but computable. It is easy to see that
    \begin{align*}
        \big\| \prod^{w}_{i=1}\exp\big(A_{p, \sigma(\check{\tau}_{i})} (\check{\tau}_{k}^{i+1}-\check{\tau}_{k}^{i})\big) \big\|
        \leq  \exp\big(\max_{q\in \sM}\|A_{p, q}\|\tau\big).
    \end{align*}
The above inequality can be used to choose the \textit{worst} switching times.
\end{rem}

\begin{rem}
	Matrix $\tilde{S}_{k}$ is not uniquely determined from \eqref{equ:xtk2}, because $\{A_{p,q}, p,q\in\sM\}$ may be not commute.
	Nevertheless, one can  guarantee that $x_{k} \in \sS_{k}$ holds for any $k\in \sN_{+}$ because inequality \eqref{eq:hold4} holds.
	Matrix  $\check{S}_{k}$ is computable, because the modes are finite, and $\tau$ is bounded. \hl{ It it easy to verify  that
		$\big\| \prod^{w}_{i=1}\exp\big(A_{p, \sigma(\check{\tau}_{i})} (\check{\tau}_{k}^{i+1}-\check{\tau}_{k}^{i})\big) \big\| \leq \exp\big(\max_{q\in \sM}\|A_{p, q}\|\tau\big)$.}
	$\tilde{x}_{t_{k+1}}$ is used to estimate the center of quantization and is dependent on the mean value of the sojourn time of the mode, which can be seen as the ``\textit{expected value}''. $\check{x}(t_{k+1})$ is used to estimate the range of quantization, which can be seen as the ``\textit{worst value}''.
\end{rem}

\subsection{Stability Analysis of  the MJLS}
In the subsection, the sufficient condition will be given to ensure the stabilization of the MJLSs under the above communication and control protocol.

For  convenience, let
\begin{align*}
	\nu_{p}        & = \max \Big\{1 - \alpha_{1, p}, \sfrac{\beta_{1, p}}{\rho_{p}}+\sfrac{\Lambda_{p}^{2}}{N^{2}}\Big\},                                                                     \\
	\alpha_{1, p}  & =\sfrac{\underline{\lambda}(Q_{p})}{\bar{\lambda}({P_{p})}}
	-\sfrac{\alpha_{p}\bar{\lambda}(S_{p}^{\top}Q_{p}S_{p})}{\underline{\lambda}({P_{p})}},                                                                                                   \\
	\beta_{1, p}   & =\Big(1+\sfrac{1}{\alpha_{p}}\Big)n\bar{\lambda}(S_{p}^{\top}Q_{p}S_{p})\Big(\sfrac{N-1}{N}\Big)^{2},                                                                    \\
	\upsilon_{p}   & = \max \Big\{\alpha_{2, p}, \sfrac{\beta_{2, p}}{\rho_{p}}+\sfrac{\Lambda_{p}^{2}}{N^{2}}\Big\},                                                                         \\
	\alpha_{2, p}  & =(1+\beta_{p})\sfrac{\bar{\lambda}(S_{p}^{\top}P_{p}S_{p})}{\underline{\lambda}({P_{p})}},                                                                               \\
	\beta_{2, p}   & =\Big(1+\sfrac{1}{\beta_{p}}\Big)n\bar{\lambda}(S_{p}^{\top}P_{p}S_{p})\Big(\sfrac{N-1}{N}\Big)^{2},                                                                     \\
	\mu_{p}        & = \max_{q\in\sM}\mu_{pq},                                                                                                                                                \\
	\mu_{pq}       & =\max\Big\{ \alpha_{3, pq}, \sfrac{\beta_{3, pq}}{\rho_{p}} \Big\},                                                                                                      \\
	\alpha_{3, pq} & = 2(1+\beta_{pq})\sfrac{\bar{\lambda}(\tilde{S}_{p}^{\top}P_{q}\tilde{S}_{p})}{\underline{\lambda}({P_{p})}}
	+\rho_{q} \chi_{p}^{2}(1+\alpha_{pq})\sfrac{1}{\underline{\lambda}(P_{p})},                                                                                                               \\
	\beta_{3, pq}  & = 2\Big(1+\sfrac{1}{\beta_{pq}}\Big)n\bar{\lambda}(\tilde{S}_{p}^{\top}P_{q}\tilde{S}_{p})\Big(\sfrac{N-1}{N}\Big)^{2} + \rho_{q}\psi_{p}^{2}(1+\sfrac{1}{\alpha_{pq}}),
\end{align*}
where  $\rho_{p}$, $\alpha_{p} $, $\beta_{p}$, $\alpha_{pq} $ and $\beta_{pq}$ are positive constants,  and $P_{p}$  are positive definite matrices, which will  be defined.
We arrive at the following result.

\begin{thm}\label{them1}
	Consider Markov jump linear system \eqref{sys:init}. Suppose that Assumptions \ref{asm:markproc} and \ref{asm:sampletime} hold. If the inequality
	\begin{equation}\label{cond:them1}
		\sum_{p\in\sM} \pi_p {\pp_{\mu, p}}\ln \mu_{p}
		+ \sum_{p\in\sM_{s}} \pi_p \pp_{\nu, p}\ln \nu_p
		+ \sum_{p\in\sM_{u}} \pi_p \pp_{\upsilon, p}\ln \upsilon_p < 0
	\end{equation}
	holds, where $\pp_{\nu, p} =e^{-2\gamma_{\hl{p}}\tau}$,
	$\pp_{\upsilon, p} =e^{-\gamma_{\hl{p}}\tau}$, and
	$\pp_{\mu, p} = 1 - e^{-2\gamma_{\hl{p}}\tau}$, then system \eqref{sys:init} reaches almost sure exponential stabilization under the communication and control protocol described in Section \ref{sec:comcont}.
\end{thm}

\begin{pf}
	Define the Lyapunov function
	\begin{equation*}
		V_p(x^{\ast}_k, E_k)=(x^{\ast}_k)^{\top}P_px^{\ast}_k + \rho_pE^2_k, \quad p \in \sM
	\end{equation*}
	which depends on the active mode $p$ ($ p=\sigma({t_{k}})$), where $P_{p}$ is a positive definite matrix, and $\rho_{p}>0$.
	By the definition of the quantization in Section \ref{sec:comcont}, it is obvious that the sequences $\{x^{\ast}_{k}, k \in \sN\}$ and $\{E_{k}, k \in \sN\}$ can be used to characterize the stability of $x(t)$.
	The proof is divided into $5$ steps as follows:

	\textit{Step 1: Sampling interval without switches occurring.}
	Two scenarios are needed to be considered as follows:

	~~\textit{a}) \textit{$(A_p, B_p)$ is stabilizable, i.e., $p \in \sM_{s}$.}
	Let $S_p = e^{(A_p + B_pK_p)\tau}$, there exists $P_p >0$ and $Q_p >0$ such that $S^{\top}_p P_p S_p - P_p = -Q_p < 0$.
	One can obtain that
	\begin{align*}
		     & (x_{k+1}^{*})^{\top} P_{p} x_{k+1}^{*}-(x_{k}^{*})^{\top} P_{p} x_{k}^{*}                                                                         \\
		=    & (S_{p} x_{k}^{*}+S_{p} \Delta_{k})^{\top} P_{p}(S_{p} x_{k}^{*}+S_{p} \Delta_{k})-(x_{k}^{*})^{\top} P_{p} x_{k}^{*}                              \\
		\leq & -(x_{k}^{*})^{\top} Q_{p} x_{k}^{*}+ \alpha  (x_{k}^{*})^{\top} S_{p}^{\top} P_{p} S_{p} x_{k}^{*}   \\
		     & +
		\sfrac{1}{\alpha}\Delta_{k}^{\top} S_{p}^{\top} P_{p} S_{p} \Delta_{k}+
		\Delta_{k}^{\top} S_{p}^{\top} P_{p} S_{p} \Delta_{k}
	\end{align*} for any $\alpha_{p} >0$,
	because $P_{p}$ is positive definite, which has the Cholesky factorization, and $2x^{\top} P_{p} y \leq x^{\top}P_{p}x + y^{\top}P_{p}y$ for any $x,y\in\sR^{n}$ and $\alpha > 0$, 
	\begin{align}\label{eq:vvp}
		     & V_p(x^{\ast}_{k+1}, E_{k+1}) \notag                                                                                       \\
		\leq & \Big(1 - \sfrac{\underline{\lambda}(Q_{p})}{\bar{\lambda}({P_{p})}}
		+\sfrac{\alpha\bar{\lambda}(S_{p}^{\top}Q_{p}S_{p})}{\underline{\lambda}({P_{p})}}\Big)(x_{k}^{*})^{\top} P_{p} x_{k}^{*} \notag \\
		     & +\Big(1+\sfrac{1}{\alpha}\Big)n\bar{\lambda}(S_{p}^{\top}Q_{p}S_{p})\Big(\sfrac{N-1}{N}\Big)^{2} E_{k}^{2}
		+\sfrac{\Lambda_{p}^{2}}{N^{2}}\rho_{p}E_{k}^{2} \notag                                                                          \\
		\leq & \nu_p V_p(x^{\ast}_k, E_k),                                                                                               
	\end{align}
	One can choose $\rho_{p}$ and $\alpha_{p}$ such that $\nu_{p}<1$ for any $p \in \sM_{s}$.

	~~\textit{b}) \textit{$(A_p, B_p)$ is unstabilizable, i.e., $p \in \sM_{u}$.}
	Let $K_{p}=\mO$ and $S_p = \exp(A_p \tau)$, one can obtain that
	\begin{equation*}
		\begin{aligned}
			  & (x_{k+1}^{*})^{\top} P_{p}x_{k+1}^{*} \\
			= & (x_{k}^{*})^{\top} S_{p}^{\top} P_{p} S_{p} x_{k}^{*}+2\left(x_{k}^{*}\right)^{\top} S_{p}^{\top} P_{p} S_{p} \Delta_{k}+\Delta_{k}^{\top} S_{p}^{\top} P_{p} S_{p} \Delta_{k}
		\end{aligned}
	\end{equation*}
	Processing similarly to \eqref{eq:vvp}, one can obtain that
	\begin{align}\label{eq:vup}
		     & V_p(x^{\ast}_{k+1}, E_{k+1}) \notag   \\
		\leq & \Big((1+\beta_{p})\sfrac{\bar{\lambda}(S_{p}^{\top}P_{p}S_{p})}{\underline{\lambda}({P_{p})}}\Big)(x_{k}^{*})^{\top} P_{p} x_{k}^{*} \notag \\
		     & + \Big(1+\sfrac{1}{\beta_{p}}\Big)n\bar{\lambda}(S_{p}^{\top}P_{p}S_{p})\Big(\sfrac{N-1}{N}\Big)^{2}E_{k}^{2}
		+\sfrac{\Lambda_{p}^{2}}{N^{2}}\rho_{p}E_{k}^{2} \notag                                                                                            \\
		\leq & \upsilon_p V_p(x^{\ast}_k, E_k),
	\end{align}
	for any $\beta_{p} > 0$.
	Obviously, $\upsilon_p > 1$ for any $\rho_{p}$ and $\beta_{p}$, $p \in \sM_{u}$

	\textit{Step 2: Sampling interval with switches occurring.}
	Notice that \eqref{eq:updaxE2}, one can obtain that
	\begin{align*}
		     & V_q(x^{\ast}_{k+1}, E_{k+1}) = \left(x_{k+1}^{*}\right)^{\top} P_{q} x_{k+1}^{*} + \rho_{q}E_{k+1}^{2}                                      \\
		\leq & (\tilde{S}_{p} [(x_{k}^{*})^{\top},(x_{k}^{*})^{\top}]^{\top}+\tilde{S}_{p} [\Delta_{k}^{\top},\Delta_{k}^{\top}]^{\top})^{\top}            \\
		     & \quad\times P_{q} (\tilde{S}_{p} [(x_{k}^{*})^{\top},(x_{k}^{*})^{\top}]^{\top}+\tilde{S}_{p} [\Delta_{k}^{\top},\Delta_{k}^{\top}]^{\top}) \\
		     & + \rho_{q} (\chi_{p} \norm{x^{\ast}_k} + \psi_{p} E_k)^{2}
	\end{align*}
	Processing similarly to \eqref{eq:vvp}, one can obtain that
	\begin{align}\label{eq:vpq}
		     & V_q(x^{\ast}_{k+1}, E_{k+1}) \notag                                                                                                                                              \\
		\leq & \Big(2(1+\beta_{pq})\sfrac{\bar{\lambda}(\tilde{S}_{p}^{\top}P_{q}\tilde{S}_{p})}{\underline{\lambda}({P_{p})}}\Big)(x_{k}^{*})^{\top} P_{p} x_{k}^{*} \notag                    \\
		     & + 2\Big(1+\sfrac{1}{\beta_{pq}}\Big)n\bar{\lambda}(\tilde{S}_{p}^{\top}P_{q}\tilde{S}_{p})\Big(\sfrac{N-1}{N}\Big)^{2}E_{k}^{2}\notag                                            \\
		     & +\rho_{q}\Big(\chi_{p}^{2}(1+\alpha_{pq})\sfrac{1}{\underline{\lambda}(P_{p})} (x_{k}^{*})^{\top} P_{p} x_{k}^{*}  + \psi_{p}^{2}(1+\sfrac{1}{\alpha_{pq}})E_{k}^{2}\Big) \notag \\
		   \leq & \mu_{pq} V_p(x^{\ast}_k,E_k), 
	\end{align}
	for any  $\alpha_{pq}>0$ and $\beta_{pq}>0$.


	\textit{Step 3: Combined bound at sampling times.}
	The sampling intervals divide the above three types. And combining them, from \eqref{eq:vvp} \eqref{eq:vup} and \eqref{eq:vpq}, one can get that
	\begin{multline}\label{eq:xkEk}
		V_{\sigma(t_k)}(x^{\ast}_{k}, E_{k}) \leq  \prod_{p, q\in\sM} \mu_{pq}^{N_{pq}} \prod_{p\in\sM_{s}} \nu_p^{N_p} \prod_{p\in\sM_{u}}\upsilon_p^{N_p} \\ \times V_{\sigma(t_0)}(x^{\ast}_0, E_0),
	\end{multline}
	where
	$N_{pq}$ denotes the number of the sampling intervals with switches occurring and $\sigma(t_{k})=p$ and $\sigma(t_{k+1})=q$,
	and $N_{p}$ denotes the number of the sampling intervals without switches occurring and $\sigma(t_{k})=p$.

	\textit{Step 4: State bound in sampling intervals.}
	Consider both switches occurring and no switch occurring scenarios in an interval $[t_{k}, t_{k+1})$.
	When no switch occurs, one can obtain that
	\begin{equation*}
		x(t) = [\mI_{n}~ \mO_{n}] \exp\big({A_{p, p}(t-t_{k})}\big) [ x^{\top}_{k}, c^{\top}_{k}]^{\top}
	\end{equation*}
	for $t \in [t_{k}, t_{k+1})$. One can obtain that
	\begin{align}\label{equ:eq20}
		\|x(t)\| & \leq \| [\mI_{n}~ \mO_{n}] \exp\big({A_{p, p}(t-t_{k})}\big) [ x^{\top}_{k}, c^{\top}_{k}]^{\top} \|\notag \\
		         & \leq \exp\big(\|{A_{p, p}\|\tau)}\big) (\|x_{k}\|+ \|c_{k}\|).
	\end{align}
	When switches occur, let$\tau_{k}^{1}$, $\tau_{k}^{2}$, $\ldots$, $\tau_{k}^{s}$
	be the switching times of the Markov process, where $ s $ is the number of switches in $[t_{k}, t_{k+1})$.
	Let $\tau_{k}^{0}=\tau_{k}$ and $\tau_{k}^{s+1}=\tau_{k+1}$, one can obtain that
	\begin{multline*}
		x(t) = [\mI_{n}~ \mO_{n}]\exp\big(A_{p, \sigma(\tau_{k}^{j})} (t-{\tau}_{k}^{j})\big)\\
		\prod^{j-1}_{i=1}\exp\big({A_{p, \sigma(\tau_{k}^{i})}(\tau_{k}^{i+1}-{\tau}_{k}^{i}})\big)
		[ x^{\top}_{k}, c^{\top}_{k}]^{\top}
	\end{multline*}
	for any $t \in [{\tau}_{k}^{j}, {\tau}_{k}^{j+1}) \subset [t_{k}, t_{k+1})$, $j \in \sN_{[0, s]}$. Processing similarly to \eqref{equ:eq20}, one can obtain that
	\begin{align}\label{equ:eq21}
		\|x(t)\| \leq \max_{q \in \sM, {p \neq q}}\exp\big(\|{A_{p, q}\|\tau}\big) (\|x_{k}\|+ \|c_{k}\|),
	\end{align}
	for any $t \in [t_{k},t_{k-1})$.
	From \eqref{equ:eq20} and \eqref{equ:eq21}, one can obtain that
	\begin{equation}\label{eq:xtb}
		\norm{x(t)} \leq \max_{q \in \sM}\exp\big(\|{A_{p, q}\|\tau}\big) (\|x_{k}\|+ \|c_{k}\|).
	\end{equation}

	\textit{Step 5: Almost sure exponential stabilization.}
	Let $\pp_{\mathrm{switch}}$ be the probability that switches occur in a sampling interval, and $1-\pp_{\mathrm{switch}}$ be the probability that not any switch occurs.
	Let
	$\pp_{\nu, p} \coloneqq \int_{2\tau}^{+\infty} \gamma_{p}e^{-\gamma_{p}\omega} \dd \omega =e^{-2\gamma_{p}\tau}$
	be the probability of the sojourn time more than $2\tau$ of  mode $p$ for $p \in \sM_{s}$.
	Let
	$\pp_{\upsilon, p} \coloneqq \int_{\tau}^{+\infty} \gamma_{\hl{p}}e^{-\gamma_{p}\omega} \dd \omega=e^{-\gamma_{p}\tau}$
	be the probability of the sojourn time more than $\tau$ of  mode $p$ for $p \in \sM_{u}$.
	Let
	$\pp_{\mu, p} \coloneqq \int^{2\tau}_{0} \gamma_{\hl{p}}e^{-\gamma_{p}\omega} \dd \omega = 1 - e^{-2\gamma_{p}\tau}$
	be the probability of the sojourn time less than $2\tau$ of  mode $p$ for $p \in \sM$.
	Obviously, 	$\pp_{\mathrm{switch}} < \sum_{p\in\sM} \pi_p \pp_{\mu, p}$.
	By using ergodic law of large numbers, from \eqref{eq:xkEk} and  \eqref{eq:xtb}, one can have that
	\begin{multline*}
		\limsup_{t \rightarrow \infty} \sfrac{1}{t} \ln \norm{x(t)}
		\leq \big (\sum_{p\in\sM} \pi_p {\pp_{\mu, p}}\ln \mu_{p}
		\\+ \sum_{p\in\sM_{s}} \pi_p \pp_{\nu, p}\ln \nu_p
		+ \sum_{p\in\sM_{u}} \pi_p \pp_{\upsilon, p}\ln \upsilon_p \big)
		k(\norm{x^{\ast}_0} + E_0)
	\end{multline*}
	in the almost sure sense, where $k$ is a positive constant.
	So, if condition \eqref{cond:them1} holds, then $\mprob{\limsup_{t \rightarrow \infty} \frac{1}{t} \ln \norm{x(t)} \leq -\varepsilon} = 1$.
	This completes the proof of Theorem \ref{them1}.
\end{pf}

\begin{rem}
	The condition of the MJLS is computable, which is independent of the explicit evolution of the Markov process.
	The generator $\Gamma$ of the Markov process, which encodes all properties of the process in a single matrix, is important in the conditions of stabilization. Different generators result in the different stabilization for the same modes.
\end{rem}

\begin{rem}
   \hl{
    The concept of dwell time  and average dwell-time have become standard assumption in the study of stability and stabilizability
    of switched and hybrid systems \citep{ieeetac2018yang:ddrc, csl2021berger}.
    \cite{auto2014liberzon:findatarate}, \cite{ieeetac2018yang:ddrc} and   \cite{ieeetac2017wakaiki:swls}  assume that
   the sampling period  is  no larger than the dwell time, that is, switches actually occur less often than once per sampling period.
   The dwell-time assumption of  switching is dropped by using the sojourn time of Markov process. The assumption of stabilizability of all individual modes is not required  in this paper.
}
\end{rem}

\subsection{Extend to the semi-MJLS case}
The semi‐MJLS is more general than the MJLS in modeling some practical systems.
In the subsection, we will extent our result to the semi-MJLS case.
The switching signal $\{\sigma(t), t \in \sR_{+}\}$ is a homogeneous semi-Markov process.
The discrete-time process $\{\sigma_{k}, k \in \sN\}$ is the embedded Markov chain of $\{\sigma(t), t \in \sR_{+}\}$ with transition probability matrix $\Lambda=[\lambda_{i j}] \in \sR^{M \times M}$ defined by $\lambda_{i j}=\mprob{\sigma_{k+1}=j | \sigma_{k}=i} >0$ if $i \neq j$, $\lambda_{ii}=0$ otherwise.
The function $F_{ij}(t)$ is a  cumulative distribution function of a sojourn time in mode $i$ before moving to mode $j$ of $\{\sigma(t), t \in \sR_{+}\}$, defined by $F_{i j}(t)\coloneqq\mathbb{P}\{h_{k+1} \leq t | \sigma_{k}=i, \sigma_{k+1}=j\}$ for any $i, j \in \sM$, $t \in \sR_{+}$.
The function $f_{ij}(t)$ is the probability density function corresponding to $F_{ij}(t)$.
The semi-Markov kernel $\Theta(t)=\left[\theta_{i j}(t)\right] \in \sR^{M \times M}$ of $\{\sigma(t), t \in \sR_{+}\}$ is defined by $ \theta_{i j}(t) \coloneqq \mprob{\sigma_{k+1}=j, h_{k+1}\leq t | \sigma_{k}=i}$, for any $i, j \in \sM$, $t \in \sR_{+}$. It it easy to check that $\theta_{i j}(t)= \lambda_{i j} F_{ij}(t)$.

    In this paper, the  cumulative distribution function of the sojourn time depends on both the current and next system mode.



\begin{assum}[semi-Markov process]\label{asm:semim}
	The semi-Markov process $\{\sigma(t), t \in \sR_{+}\}$ is irreducible and aperiodic.
\end{assum}
Similarly to the Markov process, Assumption \ref{asm:semim} implies the semi-Markov process is ergodic and has a unique stationary distribution
$\pi =[\pi_1, \pi_2, \ldots, \pi_M]$
which can be calculated by
$\pi \Lambda = \pi$ and $\sum_{i=1}^M\pi_i=1$ (see, e.g., \cite{book:smp}).

\begin{thm}\label{them2}
	Consider semi-Markov jump linear system \eqref{sys:init}. Suppose that Assumptions \ref{asm:sampletime} and \ref{asm:semim} hold. If the inequality
	\begin{equation*}\label{cond:them2}
		\sum_{p\in\sM} \pi_p \pp'_{\mu, p}\ln \mu_{p}
		+ \sum_{p\in\sM_{s}} \pi_p \pp'_{\nu, p}\ln \nu_p
		+ \sum_{p\in\sM_{u}} \pi_p \pp'_{\upsilon, p}\ln \upsilon_p < 0
	\end{equation*}
	holds, where $\pp'_{\mu, p} = \sum\limits_{q\in\sM}\lambda_{p q}\int^{2\tau}_{0}f_{pq}(\omega) \dd \omega$,
	$\pp'_{\nu, p} = \sum\limits_{q\in\sM}\lambda_{p q}\int_{2\tau}^{+\infty}f_{pq}(\omega) \dd \omega$ and
	$\pp'_{\upsilon, p}= \sum\limits_{q\in\sM}\lambda_{p q}\int_{\tau}^{+\infty}f_{pq}(\omega) \dd \omega$, then system \eqref{sys:init} reaches almost sure exponential stabilization under the communication and control protocol described in Section \ref{sec:comcont}.
\end{thm}

\begin{rem}
\hl{The difference between the Markov process and the semi-Markov process is the probability density function of the sojourn times.
    The sojourn times  of semi-Markov process are random variables with any distribution.
    The proof is the same as the proof of Theorem \ref{them1} except Step 5.
    One can deal with Step 5 by computing the probability of sojourn times via the semi-Markov kernel approach, thus obtain that $f_{i}(t)=\sum_{j \in \sM} \lambda_{i j} f_{ij}(t)$) and
    $F_{i}(t)=\sum_{j \in \sM} \theta_{i j}(t)=\sum_{j \in \sM} \lambda_{i j} F_{ij}(t)$.
    The proof is omitted for brevity.
    The Markov process can be treated as a special case of the semi-Markov process, where the probability density function $f_{i}(t) = \gamma_{i}e^{-\gamma_{i} t}$ only depends on the current system mode (see, e.g., \cite{book:smp}).
    Theorem \ref{them2} is more general than Theorem \ref{them1}.}
\end{rem}

\section{Numerical Simulation}\label{sec:ns}

\subsection{Evolution algorithm}
In this subsection, the algorithm of the state evolution is given.
Let $\Delta t$ be the time step, $end$ be the ending time of the simulation. $\Delta t$ can be designed depending on the sampling period $\tau$. Algorithm \ref{alg:alg1} shows the logic of the designed protocol to compute the states of the MJLS.

\begin{algorithm}[!htbp]
	\caption{Control and state evolution }
	\label{alg:alg1}
	\begin{algorithmic}
		\STATE Initial: System matrices $A_{i}$, $B_{i}$, $K_{i}$ for $i \in \sM$, $\Gamma$, $N$, $\Delta t$ and $end$
		\STATE Output: State evolution: $x(t)$, $\hat{x}(t)$
		\STATE Compute $\tau$ by using condition \eqref{cond:them1}
		\FOR {$t = 0: \Delta t : end$}
		\STATE Update $\sigma(t)$, $x(t)$ and $\hat{x}(t)$ by using  \eqref{sys:init}, \eqref{eq:continput} and \eqref{eq:auxilinsys}
		\IF {$t = t_{k}$ \COMMENT{$t_{k} = k\tau, k \in \sN$} }
		\IF {no switch occurs in $[t_{k-1}, t_{k})$}
		\STATE Update $E_{k}$, $x^{\ast}_{k}$ by using definition \eqref{eq:updaxE1}
		\ELSE
		\STATE Update $E_{k}$, $x^{\ast}_{k}$ by using definition \eqref{eq:updaxE2}
		\ENDIF
		\STATE Update $c_{k}$ according to the quantization.
		\ENDIF
		\ENDFOR
	\end{algorithmic}
\end{algorithm}

\subsection{Numerical Examples}
In this subsection, some numerical examples are provided to demonstrate the validity of the obtained theoretical results.

Consider a MJLS with three modes with system matrices:
\begin{align*}
	A_{1} & =\begin{bmatrix}
		1 & 0  \\
		0 & -1
	\end{bmatrix},
	B_{1} = \begin{bmatrix}
		1 \\
		0
	\end{bmatrix},  K_{1}=\begin{bmatrix}
		-2 & 0
	\end{bmatrix}, \\
	A_{2} & =\begin{bmatrix}
		1 & 0  \\
		0 & -1
	\end{bmatrix},
	B_{2} = \begin{bmatrix}
		0 \\
		1
	\end{bmatrix},  K_{2}=\begin{bmatrix}
		0 & 0
	\end{bmatrix}, \\
	A_{3} & =\begin{bmatrix}
		0  & 1 \\
		-1 & 0
	\end{bmatrix},
	B_{3} = \begin{bmatrix}
		0 \\
		1
	\end{bmatrix},  K_{3}=\begin{bmatrix}
		0 & -4
	\end{bmatrix},
\end{align*}
and
generator of the Markov process
\begin{equation*}
	\Gamma = \begin{bmatrix}
		-0.050 & 0.010 & 0.040  \\
		0.075  & -0.15 & 0.075  \\
		0.035  & 0.005 & -0.045
	\end{bmatrix}.
\end{equation*}

\begin{figure*}[!htbp]
    \centering
    \subfigure[Time evolution of state $x_{1}(t)$ and $\hat{x}_{1}(t)$]{\includegraphics[width=0.47\linewidth]{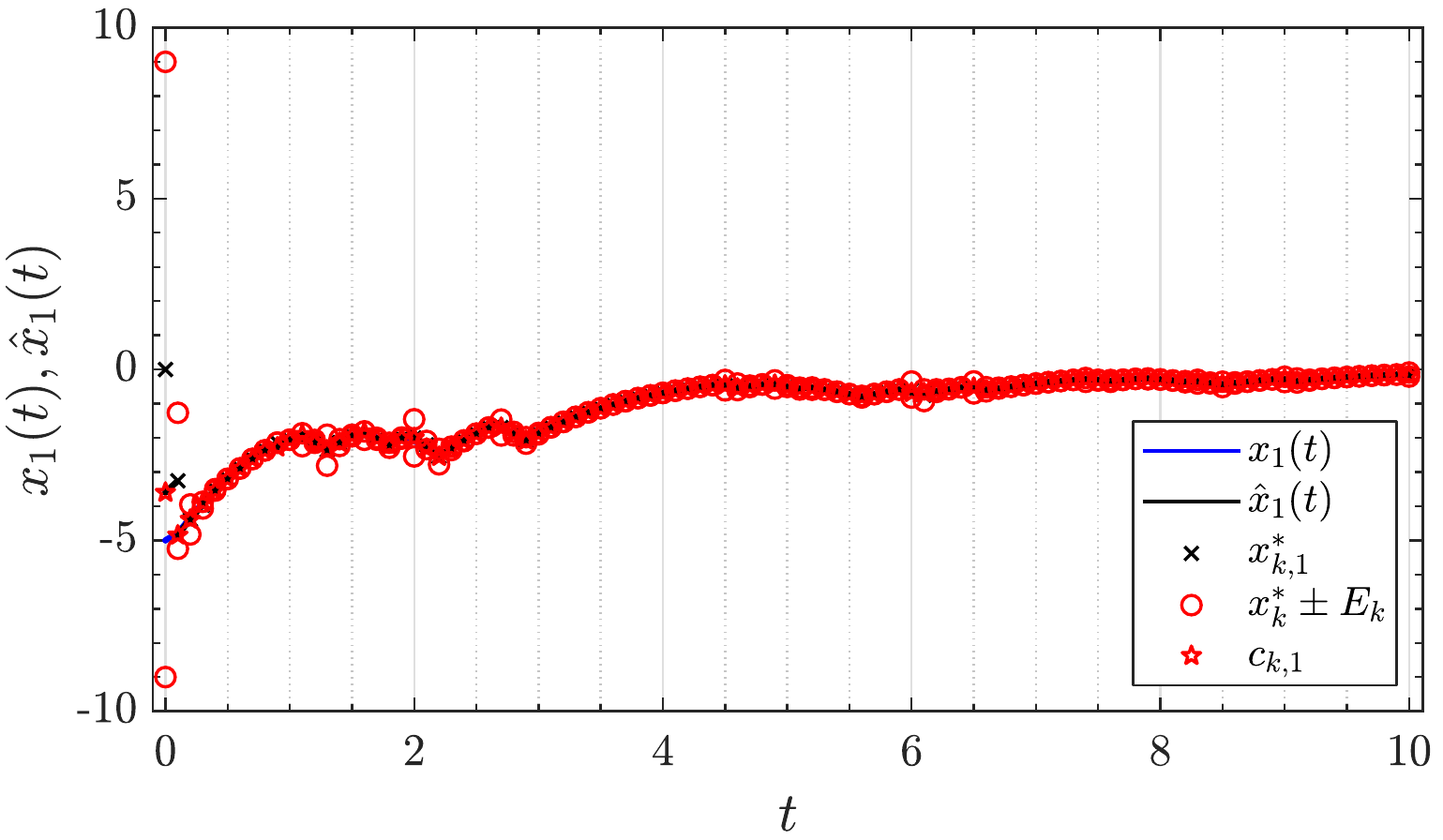}}~~~~
    \subfigure[Time evolution of state $x_{2}(t)$ and $\hat{x}_{2}(t)$]{\includegraphics[width=0.47\linewidth]{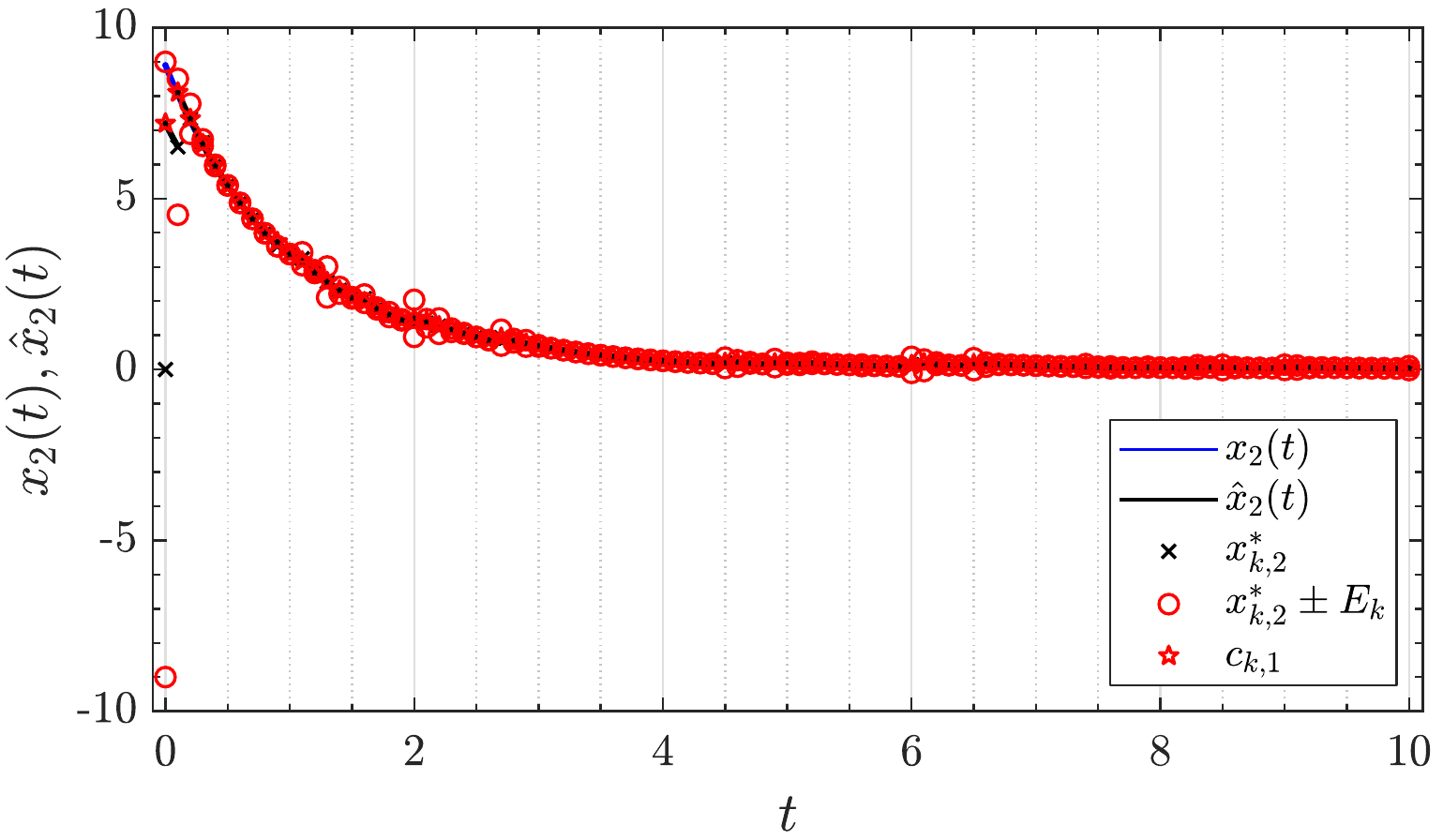}}

    \subfigure[The quantization $x_{1}(t)$]{\includegraphics[width=0.47\linewidth]{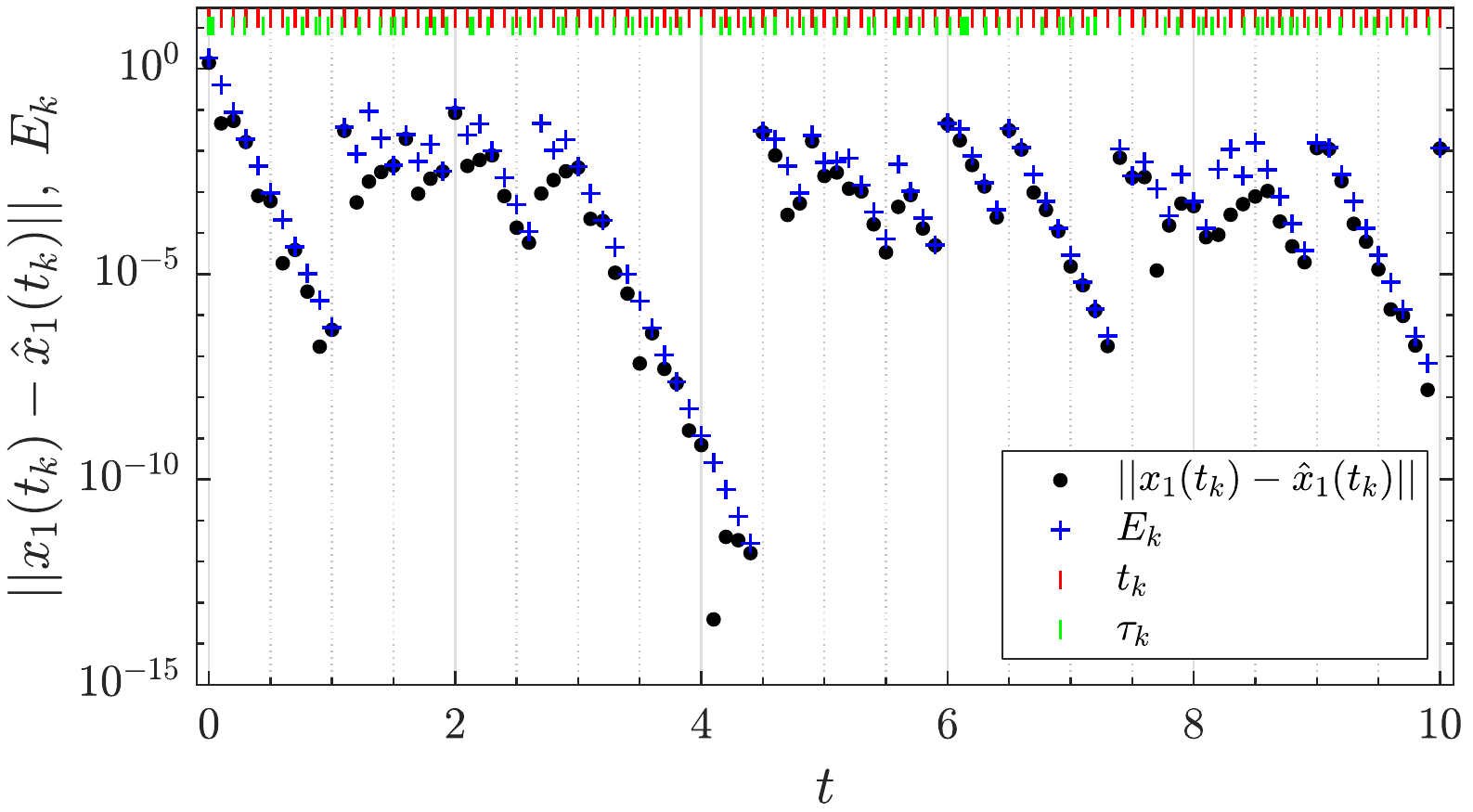}}~~~~
    \subfigure[The quantization $x_{2}(t)$]{\includegraphics[width=0.47\linewidth]{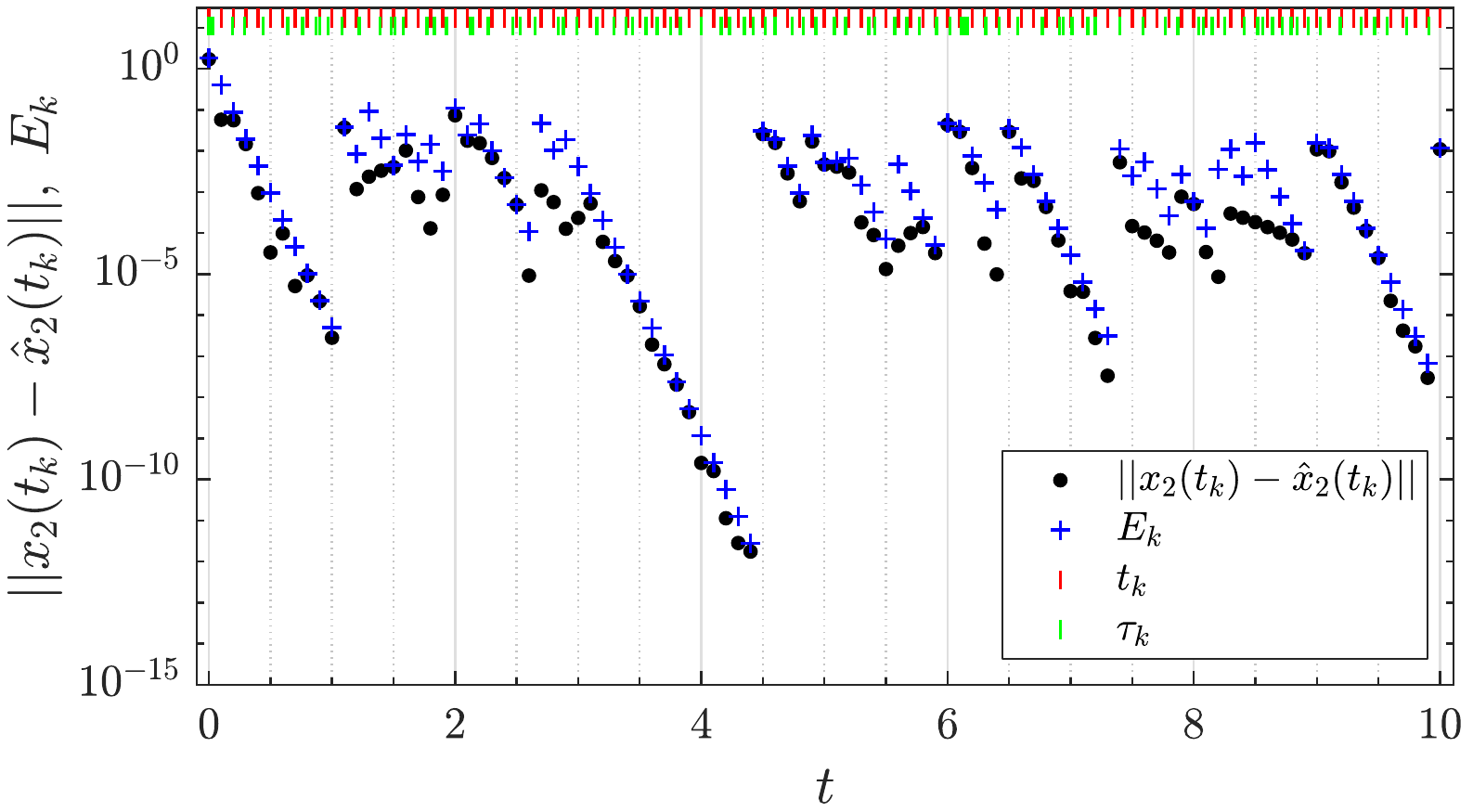}}
    \caption{The states and quantization of the MJLS with a Markov process sample path. }
    \label{fig3}
\end{figure*}

First, consider the MJLS under the protocol defined in Section \ref{sec:comcont} and \ref{sec:uq}.
Choose $P_{1} = \mI_{n}$, $P_{2} = \mI_{n}$, for $p \in \{1,2,3\}$, $\tau=0.1$, $N=10$ and $\Delta t=0.0001$ such that the conditions of Theorem \ref{them1} are satisfied (the detailed results are shown in Table \ref{table:exam1}).
The evolution of the states of the MJLS
is shown in Fig. \ref{fig3} for the chosen initial values $[-5, 8.9]^{\top}$.
One can observe from Fig. \ref{fig3} that stabilization can be reached.


\begin{table}[!htp]
	\center

	\begin{tabular}{cccccc}
		\hline
		$p$ & stabilizability & $\nu_p$ & $\upsilon_{p}$ & $\mu_p$ & $\pi_{p}$ \\
		\hline
		1   & Yes             & 0.8766  & /              & 5.5734  & 0.4386    \\
		2   & No              & /       & 1.3435         & 4.8951  & 0.1404    \\
		3   & Yes             & 0.9883  & /              & 7.7885  & 0.4211    \\
		\hline
	\end{tabular}
	\caption{The detailed computation results.}
	\label{table:exam1}
\end{table}

\begin{figure}[!htbp]
	\centering
	\includegraphics[width=0.9\linewidth]{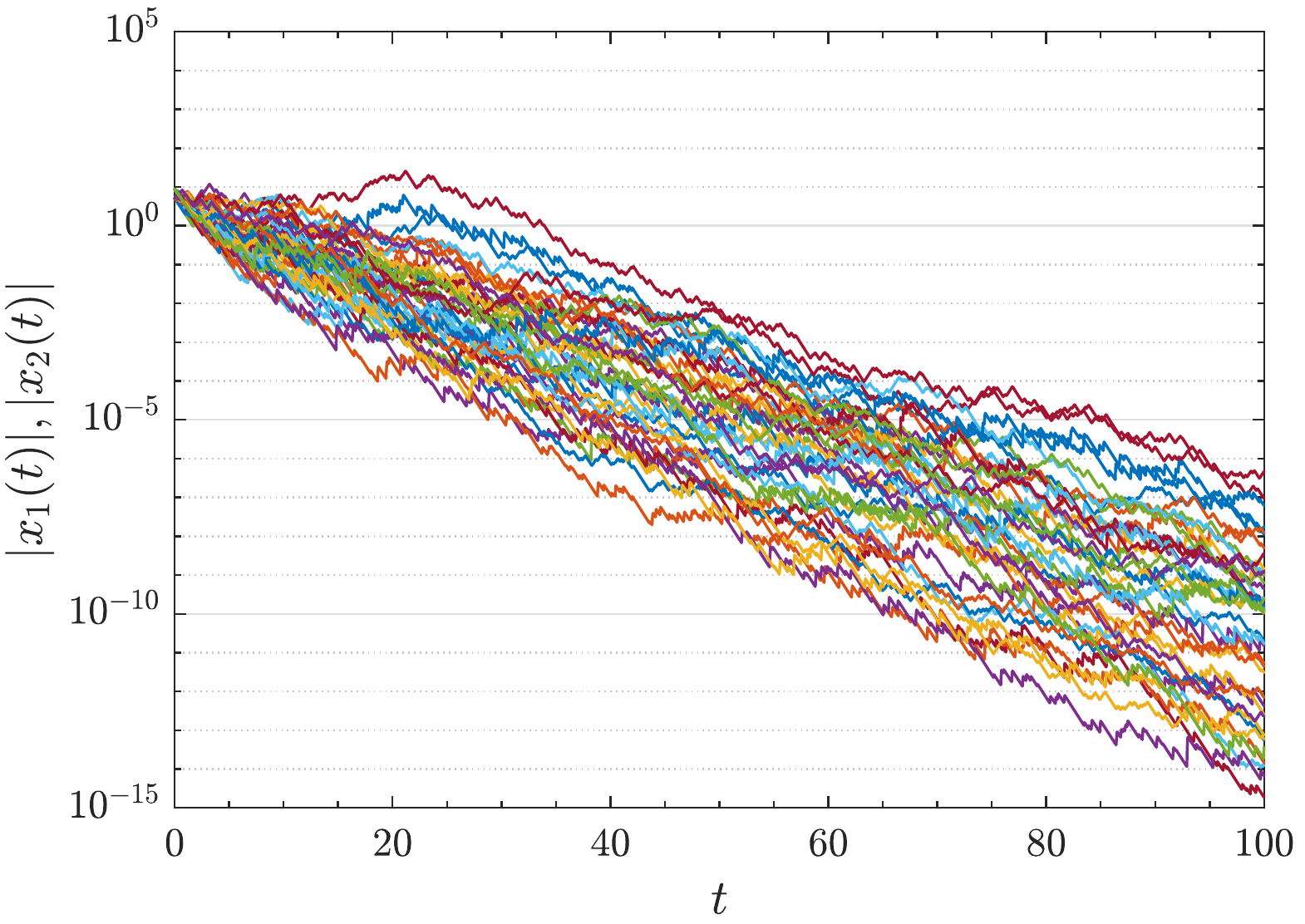}
	\caption{The states of 20 realizations with $x(0) = [-5, 8.9]^{\top}$.}
	\label{fig:real100}
\end{figure}

Next, the $20$ realizations of the MJLS  are given and Fig. \ref{fig:real100} illustrates the  $|x_{1}(t)|$ and $|x_{2}(t)|$ starting with the same initial value $[-5, 8.9]^{\top}$ and generator $\Gamma$. Apparently, the MJLS under the designed communication and control protocol achieves almost sure exponential stabilization.

\section{Conclusion}\label{sec:conclusion}

In this paper, we consider the stabilization problem of the MJLSs under the communication data-rate constraints,
where the switching signal is a continuous-time Markov process.
Sampling and quantization are used as fundamental tools to deal with the problem.
Under the proposed communication and control protocol, the sufficient conditions are given to ensure the almost sure exponential stabilization of the MJLSs.
The conditions depend on the generator of the Markov process.
The sampling times and the jump time is also independent.
We extend the result to the semi-MJLSs case.

In future, we will extend our results to Markov/semi-Markov jump nonlinear systems. The communication and control protocols will applied to networked control systems. Typical communication channels are noisy and have delays. The noise and delays will be considered in stabilization of MJLSs/semi-MJLSs.

\begin{ack}
	This work was supported
	in part by the National Natural Science Foundation of China under Grant 61873171, 61872429 and 61973241,
	in part by the Natural Science Foundation of Guangdong Province, China under Grant 2019A1515012192.
\end{ack}

\end{document}